\documentclass[letter,12pt]{amsart}
\usepackage{latexsym}
\usepackage{graphics}
\usepackage{amssymb}
\usepackage{amsmath}
\usepackage{epsfig}
\usepackage{tikz}%%
\usepackage{amscd}
\usepackage[all]{xy}

\def\sm{Smarandache\,}

\newcommand{\RR}{\mathbb{R}}
\newcommand{\RQ}{\mathbb{Q}}

\newcommand{\ZZ}{\mathbb{Z}}

\begin{document}
  \vspace{10mm}

\title[Right feeble groups]{Right feeble groups }
\bigskip
\author{Hiba F. Fayoumi and  Hee Sik Kim}
\maketitle
\bigskip
\begin{abstract}
Right feeble groups are defined as groupoids $(X,*)$ such that (i) $x, y\in X$ implies the existence of $a, b \in X$ such that $a*x = y$ and $b*y = x$. Furthermore, (ii) if $x, y, z \in X$ then there is an element $w\in X$ such that $x*(y*z) = w*z$.  These groupoids have a ``remnant" group structure, which includes many other groupoids.  In this paper, we investigate some properties of these groupoids. Enough examples are supplied to support the argument that they form a suitable class for systematic investigation.
\end{abstract}
\bigskip

\renewcommand{\thefootnote}{}


\section{ Introduction}
\medskip\medskip

Among the many generalizations of the idea of groups already in the literature \cite{AKN16, Hoo13, WK18} we believe that the class of right feeble groups has not been discussed in its undisguised form until now.  The terminology ``feeble" comes from the fact that conditions (i) and (ii) defined in section 3 are replacements of the existence of inverses and the associative law.
Indeed, suppose that the right feeble group $(X,*)$ contains an identity element $e$, where $x*e = e*x = x$ for all $x\in X$. Thus $a*x = e$ and $b*e= x$ implies $b= x$ and $a= x^{-1}$, i.e., the existence of inverses is assured. Also rule (ii) reads that for any $x, y, z$ in the groupoid $(X,*)$ there is an element $w \in X$ such that $x*(y*z) = w*z$. For groups we find that $w = x*y$ and thus rule (ii) is a more general version of the associative law.
\medskip

\section{  Preliminaries}
\medskip

Amid  the various kinds of groupoids, leftoids play a special role. For instance, let $f: X \to X$ be any function, we define $(X,*, f)$ by the formula: $x*y :=f(x)$.  Thus, if we consider the groupoid-product $(X,*, f)\Box (X,\bullet, g):= (X,\Box)$, with $x\Box y = (x*y)\bullet(y*x) = f(x)\bullet g(y) = g(f(x))$, it follows that $(X, \Box) = (X, \Box, g\circ f)$ is also a leftoid model of the composition of functions $(g\circ f)(x) = g(f(x))$.  If $f(x) = x$, then $x*y = f(x) = x$ produces the left zero semigroup for which we obtain $(x*y)\bullet(y*x) = x\bullet y = x\Box y$ if $(X,*)$ is the left zero semigroup.  If $(X,\bullet)$ is the left zero semigroup, then $(x*y)\bullet (y*x) = x*y = x\Box y$ as well. Meaning, the left zero semigroup acts like a multiplicative identity in $(Bin(X), \Box)$ (\cite{KN08}). H. Fayoumi (\cite{Fa11}) introduced the notion of the center $ZBin(X)$ in the semigroup $Bin(X)$ of all binary systems on a set $X$, and showed that if $%
(X,\bullet )\in ZBin(X)$, then $x\not=y$ implies $\{x,y\}=\{x\bullet
y,y\bullet x\}$. Moreover, she showed that a groupoid $(X,\bullet )\in ZBin(X)$
if and only if it is a locally-zero groupoid.

Suppose that $f: X \to X$ is a surjection.  If $x, y\in X$, then $f(a) = y$ and $f(b) = y$ for some choice of $a$ and $b$, i.e., $a*b = f(a) = y$ and $b*y = f(b) = x$, so that condition (i) holds. Also $x*(y*z) = f(x)$ and $w*z = f(w)$.  Since $x = w$ guarantees that $f(x) = f(w)$, it follows that condition (ii) holds as well and that the leftoid $(X,*, f)$ is a right feeble group.

Given a set $X$ and a function $f:X\to X$, we consider a
groupoid $(X,*,f)$, where $x*y := f(x)$ for any $x, y\in X$.  Such a groupoid is called a {\it leftoid} over $f$ (\cite{KN08}).  Similarly, we define a rightoid, i.e., $x*y := g(y)$ for all $x, y\in X$, where $g: X \to X$ is a function.
Another idea of interest that will be useful in what follows is that of Smarandache disjointness.  Given algebra types $(X,*)$ (type-$P_1$) and $(X,\circ)$
(type-$P_2$), we shall consider them to be {\it \sm disjoint}  (\cite{AKN05}) if
the following  two conditions hold:
\begin{itemize}
 \item[\rm (A)] If $(X,*)$ is a type-$P_1$-algebra with $|\, X\, |> 1$ then it cannot
 be a \sm-type-$P_2$-algebra $(X,\circ)$;
 \item[\rm (B)] If $(X,\circ)$ is a type-$P_2$-algebra with $|\, X\, |> 1$ then it cannot
 be a \sm-type-$P_1$-algebra $(X,*)$.
\end{itemize}
Thus, if  $\emph{K}_1$ and $\emph{K}_2$ are two classes of mathematical objects, it may be that $
\emph{K}_1\cap \emph{K}_2$ consists precisely of one single object.
Frequently this one single object is ``trivial" in some way.  For example, let $\emph{K}_1$ be the class of $d$-algebras (\cite{NK99}) and let $\emph{K}_2$ be the class of groups. If $(X,*, 0)$ is both a $d$-algebra and a group with identity $e$, then $0*x = 0$ implies $x= e$, the identity of the group, and thus $X=\{0\}$. Hence, $e= 0$ as well and $e*e = e = 0 = 0*0$.

Note that we may enlarge $\emph{K}_1$ to the class of all groupoids $(X,*,0)$ for which $0*x = 0$ for all $x\in X$, to obtain the same conclusion.  Similarly, $\emph{K}_2$ may be enlarged to the class of groupoids $(X,*, e)$ where $a*x = a$ implies $x = e$.  Hence $\emph{K}_1\cap \emph{K}_2$ consists of the single groupoid $(X=\{u\}, u*u = u)$, and $\emph{K}_1$ and $\emph{K}_2$ are then Smarandache disjoint.

If $(X,*)$ is both a leftoid and a rightoid, then $x*y = C$, a constant from $X$.  In this case, if $(X,*)$ has $x*y = C$ and $(X, \bullet)$ has $x\bullet y = D$, then the groupoids are isomorphic.  Indeed, let $\varphi: X \to X$ be any bijection such that $\varphi(C) = D$, so that $\varphi(x*y) = \varphi(C) = D = \varphi(x)\bullet \varphi(y)$.
Hence, leftoids and rightoids on a set $X$ are Smarandache disjoint up to isomorphism of groupoids.
\medskip

\section{ Right feeble groups}
\medskip\medskip

A groupoid $(X,*)$ is said to be a {\it right feeble group} if
\begin{itemize}
 \item[{\rm (i)}] for any $x, y\in X$, there exist $a, b\in X$ such that $a*x = y,$ and $ b*y = x$,
 \item[{\rm (ii)}] for any $x, y, z \in X$,  there exists $w \in X$  such that $x*(y*z) = w*z.$
\end{itemize}
\medskip

{\bf Example 3.1.} Let $\RR$ be the set of all real numbers. If we define
a binary operation ``$*$" on $\RR$ by
$x*y:= \frac{1}{2}(x + y)$ for any $x, y\in \RR$, then $(\RR,*)$ is a right feeble group.  In fact, given $x, y\in \RR $, if we take $a:= 2y - x,$ and $ b:= 2x -y$, then $a*x = y$ and $b*y = x$. Since $x*(y*z) = \frac{1}{2}(x + \frac{y+z}{2})$, we let $w:= x + \frac{1}{2}y - \frac{1}{2}z$, then $x*(y*z) = w*z$.
\medskip

Note that $(\RR, *)$ in Example 3.1 is neither a group nor a semigroup.
Assume that $(\RR, *)$ is a group with identity $e$.  Then $x*e = x$ for all $x\in \RR$. It follows that $\frac{x + e}{2} = x$ for all $x\in \RR$, which shows that $x =e$ for all $x\in \RR$, i.e., $|\RR|=1$, a contradiction.
Moreover, $1*(3*5) = \frac{5}{2}\not = \frac{7}{2} = (1*3)*5$.
\medskip

{\bf Example 3.2.}  Let $(X, +, \cdot)$ be a field and let $\alpha, \beta, \gamma \in X$. If we define a binary operation ``$*$" on $X$ by
$x*y:= \alpha + \beta x + \gamma y$ for any $x, y\in X$, then $(X,*)$ is a right feeble group.

In fact, for any $x, y\in X$,  if we let $a:= \frac{1}{\beta}(y-\alpha -\gamma x)$, and $b:= \frac{1}{\beta}(x - \alpha -\gamma y)$, then it is easy to see that  $a*x = y$ and $b *y = x$. Given $x, y, z, w\in X$, since $x*(y*z) = \alpha(1 + \gamma) +\beta(x +y) + \gamma^2z$ and $w*z = \alpha + \beta w + \gamma x$, if we take $w:= \frac{1}{\beta}[\alpha\gamma + \beta(x +y) + \gamma(\gamma -1)x]$, then $x*(y*z) = w*z$.  This shows that $(X,*)$ is a right feeble group.
\medskip

{\bf Proposition 3.3.} {\sl Every group is a right feeble group.}
\medskip

{\it Proof.} Let $(X,*)$ be a group with identity $e$.  Given $x, y\in X$, if we take $a:= y * x^{-1},$ and $ b:= x* y^{-1}$, then $a*x = (y*x^{-1})*x = y*(x^{-1}* x) = y$ and $b*y = (x* y^{-1})*y = x*e = x$.  Given $x, y, z \in X$, if we let $w:= x*y$, then $x*(y*z) = (x*y)*z = w*z$. Once again proving that $(X,*)$ is a right feeble group. \qed
\medskip

{\bf Proposition 3.4.} { Let $(X,*)$ be a leftoid for $\varphi$. If $\varphi(X) = X$, then $(X,*)$ is a right feeble group.}
\medskip

{\it Proof.} Given $x, y\in X$, since $\varphi$ is onto, there exists $a\in X$ such that $\varphi(a) = y$. Since $(X,*)$ is a leftoid for $\varphi$, we have $a*x = \varphi(a) = y$.  Similarly, if we take $b\in X$ such that $\varphi(b) = x$, then $b*y = \varphi(b) = x$.

Given $x, y, z\in X$, since $(X,*)$ is a leftoid for $\varphi$, we have $x*(y*z) = \varphi(x) = x * z$. \qed
\medskip

The notion of Smarandache was introduced by Smarandache and Kandasamy in (\cite{Kan02}) studied the concept of Smarandache groupoids and Smarandache Bol groupoids. Padilla in (\cite{Pad88})
examined Smarandache algebraic structures. Allen, Kim and Neggers (\cite{AKN05}) introduced
Smarandache disjointness in $BCK/d$-algebras. For more information on the notion of Smarandache we refer to (\cite{Kan02}).\\
In the next theorem, we consider the class of groups and the class of leftoids.

\medskip

{\bf Theorem 3.5.} {\sl The class of groups and the class of leftoids are Smarandache disjoint.}
\medskip

{\it Proof.} Let $(X,*)$ be both a leftoid for $\varphi$ and a group.
Then $e =x*x^{-1} = \varphi(x)$ for any $x\in X$, where $e$ is the identity for the group $(X,*)$. It follows that $x*y = \varphi(x) = e = x*x$ for any $x, y\in X$.  Since $(X,*)$ is a group, we obtain $x = y$ for all $x, y\in X$, proving that $|X|=1$. \qed
\medskip

{\bf Proposition 3.6.} {\sl Let $(X,*)$ be a right feeble group. If $f: (X,*) \to (Y, \bullet)$ is an epimorphism of groupoids, then $(Y, \bullet)$ is a right feeble group.}
\medskip

{\it Proof.} Given $x, y\in Y$, since $f$ is onto, there exist $a, b\in X$ such that $x = f(a),$ and $ y = f(b)$. Since $(X,*)$ is a right feeble group, there exist $p, q, r, s \in X$ such that $p*a = b, q*b = a, r*b= a,$ and $ s*a = b$.  It follows that $y = f(b) = f(p*a) = f(p)\bullet f(a)= f(p)\bullet x$ and $x = f(a) = f(q*b) = f(q)\bullet f(b) = f(q) \bullet y$.

Given $x, y, z\in Y$, since $f$ is onto, there exist $a, b, c\in X$ such that $x= f(a), y = f(b),$ and $ z = f(z)$.  Since $(X,*)$ is a right feeble group, there exists $w\in X$ such that $a*(b*c) = w*c$.  It follows that $x\bullet(y\bullet z) = f(a) \bullet (f(b) \bullet f(c))= f(a*(b*c)) = f(w*c) = f(w)\bullet f(c) = f(w) \bullet z$, proving that $(Y, \bullet)$ is also a right feeble group. \qed
\medskip

{\bf Proposition 3.7.} {\sl Let $(X,*),$ and $ (Y, \bullet)$ be right feeble  groups and let $Z:=X\times Y$. Define $(x, y)\triangledown (u, v):= (x*u, y\bullet v)$ for all $(x, y), (u, v)\in Z$.  Then $(Z, \triangledown)$ is also a right feeble group.}
\medskip

{\it Proof.} Straightforward. \qed
\medskip

{\bf Example 3.8.} Let $\RR$ be the set of all real numbers and ``$+$" be the usual addition on $X$.  Then $(\RR,+)$ forms a group. By Proposition 3.2, it is a right feeble group.  Let $A:=[0, \infty)$.  Then $(A,+)$ is a subgroupoid of $(\RR,+)$, but not a right feeble group.  In fact, if we assume $(A,+)$ is a right feeble group, then, for any $x, y\in A$, there exist $a, b\in \RR$ such that $a+x = y,$ and $ b + y = x$. It follows that $y= a + x = a + (b +y)$. Since $a, b\in [0, \infty)$, we obtain
$a = b = 0$, proving that $x = b +y = 0 +y =y$ for all $x, y\in A$, a contradiction. \qed
\medskip

In Example 3.8, the subgroupoid $(A,+)$ is not a right feeble group.
Hence, the following question:  If $(A,*)$ is a subgroupoid of a right feeble group $(X,*)$, under what condition(s) will $(A,*)$ be a right feeble group?  \\
To solve this problem, we introduce the notion of ``divisibility".
\medskip

Let $(X,*)$ be a groupoid. A subgroupoid $(A,*)$ is said to be {\it divisible} in $(X,*)$ if $a*x = y$ and $x, y\in A$, then $a\in A$.
\medskip

{\bf Example 3.9.} (a) Let $\RR$ be the set of all real numbers and let $\RQ$ be the set of all rational numbers without 0.  Define a binary operation ``$*$" on $\RR$ by $x*y:= xy$ (the usual multiplication).  Assume that $a*x = y$
and $x, y\in \RQ$.  Since $x\not = 0$, we have $a = yx^{-1}\in \RQ$, which shows that $(\RQ,*)$ is divisible.
\smallskip

\noindent
(b) Let $\RR$ be the set of all real numbers and let $\ZZ$ be the set of all integers.  Then $\frac{1}{2}\cdot 4 = 2$ and $4, 2 \in \ZZ$, but $\frac{1}{2} \not \in \ZZ$, which shows that $(\ZZ,*)$ is not divisible.
\medskip

If $(A,*)$ and $(B,*)$ are divisible in a groupoid $(X,*)$, and if $a*x = y$, where $x, y\in A\cap B$, then $a\in A\cap B$, and $A\cap B\not = \emptyset$ implies $(A\cap B, *)$ is divisible in $(X,*)$.
\medskip

{\bf Theorem 3.10.} {\sl Let $(X,*)$ be a right feeble group. If $(A, *)$ is divisible in $(X,*)$, then $(A,*)$ is a right feeble group. }
\medskip

{\it Proof.} (i) Given $x, y\in A$, since $(X,*)$ is a right feeble group, there exist $a, b \in X$ such that $a*x = y,$ and $b*y = x$.  Since $(A,*)$ is divisible and $x, y\in A$, we have $a, b \in A$ such that $a*x = y,$ and $b*y = x$.

\noindent
(ii) Given $x, y, z\in A$, we let $u:= x*(y*z)$.  Since $(A,*)$ is a subgroupoid of $(X,*)$, we have $u\in A$. Since $A$ is divisible and $z\in A$, there exists $w\in A$ such that $w*z = u = x*(y*z)$. Hence, $(A,*)$ is a right feeble group. \qed
\medskip

\section{ Right entire and right asymmetric}
\medskip

Let $(X, *)$ be a groupoid. We define a set $\rho(X,*)$ by
\[
\rho(X,*):= \{x\in X| X*x = X\}.
\]

{\bf Proposition 4.1.} {\sl If $(X,*)$ is a right feeble group, then
$\rho(X,*)= X$.}
\medskip

{\it Proof.} For any $x, y\in X$, since $(X,*)$ is a right feeble group, there exist $a, b \in X$ such that $y = a*x, x = b*y$. It follows that $y = a*x \in X*x$. Hence $X\subseteq X*x$, i.e., $X= X*x$, for all $x\in X$.
This means that $x\in \rho(X,*)$ for all $x\in X$, proving that $X=\rho(X,*)$. \qed
\medskip

The groupoid $(X,*)$ discussed in Proposition 4.1 is said to be {\it right entire}.  It follows immediately that right feeble groups are right entire groups.
\medskip

{\bf Proposition 4.2.}  {\sl Let $(X,*), (Y, \bullet)$ be right entire  groups and let $Z:=X\times Y$. Define $(x, y)\triangledown (u, v):= (x*u, y\bullet v)$ for all $(x, y), (u, v)\in Z$.  Then $(Z, \triangledown)$ is also a right entire groupoid.}
\medskip

{\it Proof.} Straightforward. \qed
\medskip

{\bf Proposition 4.3.} {\sl Let $(X,*)$ be a right entire groupoid.  If $f: (X, *)\to (Y, \bullet)$ is an epimorphism of groupoids, then $(Y, \bullet)$ is a right entire groupoid.}
\medskip

{\it Proof.} Given $x, y\in Y$, since $f$ is onto, there exist $p, q\in X$ such that $x = f(q), y = f(p)$.  Since $(X,*)$ is right entire, there exists $r\in X$ such that $r*p = q$ and hence $x = f(q) = f(r*p) = f(r)\bullet f(p) = f(r)\bullet y$.  This shows that $Y= Y\bullet y$ for all $y\in Y$, proving that $\rho(Y, \bullet) = Y$. \qed
\medskip

{\bf Proposition 4.4.} {\sl If $(X,*)$ is a leftoid for $\varphi$ and right entire, then $\varphi(X) = X$.}
\medskip

{\it Proof.} Let $(X,*)$ be a leftoid for $\varphi$ and a right entire groupoid.  Then $X*x = X$ for all $x\in X$, i.e., there exists $b\in X$ such that $a = b*x$ for any $a\in X$.  Since $(X,*)$ is a leftoid for $\varphi$, we have $a= b*x = \varphi(b) \in \varphi(X)$ for all $a\in X$, proving that $X\subseteq \varphi(X)$. \qed
\medskip

A groupoid $(X,*)$ is said to be {\it right asymmetric} if $a*x = y,$ and $b*y = x$ for some $a, b \in X$, then $x =y$.
\medskip

{\bf Example 4.5.} Let $X:=[0, \infty)$ and let $x*y:= x + y$ for all $x, y \in X$. Then $(X,*)$ is right asymmetric. In fact, if $a*x =y,$ and $b*y = x$, then $a + x = y,$ and $b +y = x$, and hence $y = a + x = a + (b + y) = (a + b) + y$. It follows that $a + b = 0$, i.e., $a = b = 0$. This shows that $x=y$.
\medskip

{\bf Theorem 4.6. } {\sl The class of right entire groupoids and the class of right asymmetric groupoids are Smarandache disjoint.}
\medskip

{\it Proof.} Assume $(X,*)$ be both a right entire groupoid and a right asymmetric groupoid.  Then $X*x = X, X*y = X$ for all $x, y\in X$. It follows that $a*x = y,$ and $b*y = x$ for some $a, b\in X$.  Since $(X,*)$ is right asymmetric, we obtain $x=y$ for all $x, y\in X$.  Hence, $|X|=1$.
\qed
\medskip

\section{Some relations }
\medskip

Given a groupoid $(X,*)$, we define a binary operation ``$\leq$" on $X$ by
\[
x\leq y\,\, \iff\,\, \exists\, a\, \in X\,\, {\rm s.t.} \,\, a*x = y.
\]
\medskip

{\bf Proposition 5.1.} {\sl If $(X,*)$ is a right entire groupoid, then $\leq$ is reflexive.}
\medskip

{\it Proof.} Assume $(X,*)$ is right entire.  Then $\rho(X,*) = X$, i.e., $X*x = X$ for all $x\in X$.  It follows that there exists $a\in X$ such that $x = a*x$  for any $x\in X$, i.e., $x\leq x$ for any $x\in X$. \qed
\medskip

{\bf Proposition 5.2.} {\sl A groupoid $(X,*)$ is  right asymmetric if and only if
 $\leq$ is anti-symmetric.}
\medskip

{\it Proof.} Let $(X,*)$ be a right asymmetric groupoid.  Assume $x\leq y$ and $y\leq x$.  Then there exist $a, b\in X$ such that $a*x = y,$ and $ b*y = x$.
Since $(X,*)$ is right asymmetric, we obtain $x = y$.
The converse is trivial and we omit the proof. \qed
\medskip

{\bf Proposition 5.3.} {\sl If $(X,*)$ is a right feeble group, then $\leq$ is transitive.}
\medskip

{\it Proof.} Assume that $x\leq y,$ and $y\leq z$. Then there exist $a, b\in X$ such that $a*x = y,$ and $b*y = z$.  Since $(X,*)$ is a right feeble group, there exists $c\in X$ such that $z = b*y = b*(a*x) = c*x$, which shows that $x\leq z$.  \qed
\medskip

{\bf Proposition 5.4.} {\sl Let $(X,*)$ be a groupoid with $e\in X$ such that $e*x = x$ for all $x\in X$.  Then $x\leq x$ for all $x\in X$.}
\medskip

Note that  `$(X,*)$ has a left identity' does not mean $X*x = X$ for some $x\in X$.  For example, let $\emph{N}:=\{0, 1, 2, \cdots\}$.  Then
$(\emph{N}, +)$ has an identity 0 and  $0+x = x$ for all $x\in \emph{N}$, but $\emph{N} + 2 =\{2, 3, \cdots \} \not = \emph{N}$.
\medskip

{\bf Proposition 5.5.} {\sl Let $(X,*)$ be a groupoid and let $a\in X$ such that $a*X = X$. Then there exists $x\in X$ such that $x\leq y$ for any $y\in X$.}
\medskip

{\it Proof.} If $a*X = X$, then there exists $x \in X$ such that $y= a*x$ for any $y\in X$. It follows that $x\leq y$. \qed
\medskip

Let $\emph{N}:=\{0, 1, 2, \cdots\}$. Then $0+ \emph{N} = \emph{N}$ and hence $y\leq y$ for all $y\in \emph{N}$.
\medskip

{\bf Proposition 5.6.} {\sl Let $(A,*)$ be a divisible subgroupoid of a groupoid $(X,*)$. If $x, y\in A$ such that $x\leq y$ in $(X,*)$, then $x\leq y$ in $(A,*)$.}
\medskip

{\it Proof.}  Let $x, y\in A$ such that $x\leq y$ in $(X,*)$. Then there exists $a\in X$ such that $a*x = y$.  Since $(A,*)$ is divisible and $x, y\in A$, we have $a\in A$, i.e., $x\leq y$ in $(A,*)$. \qed
\medskip

{\bf Theorem 5.7.} {\sl Let $(X,*)$ be a right entire groupoid.  If $(A,*)$ is a divisible subgroupoid of $(X,*)$, then $(A,*)$ is right entire.}
\medskip

{\it Proof.} Given $a\in A$, since $(X,*)$ is right entire, we have $X*a = X$. It follows that there exists $y\in X$ such that $x = y*a$ for any $x\in A$. Since $(A,*)$ is divisible, we obtain $y\in A$ and hence $x = y*a \in A * a$.  Hence $A\subseteq A*a$ for all $a\in A$. Clearly, $A*a \subseteq A$, proving that $A= A*a$ for any $a\in A$. \qed
\medskip

{\bf Proposition 5.8.} {\sl Let $(X,*)$ be a right asymmetric groupoid.  If $(A,*)$ is a divisible groupoid in $(X,*)$, then $(A,*)$ is also right asymmetric.}
\medskip

{\it Proof.} If $x, y\in A$, then $x, y\in X$.  Since $(X,*)$ is right asymmetric, if $a*x = y,$ and $b*y = x$ for some $a, b\in X$, then $x=y$.
We show that $a, b\in A$. Consider $a*x = y$.  Since $x, y\in A$ and $A$ is divisible, we obtain $a\in A$.  Similarly, $b\in A$. Hence $(A,*)$ is right asymmetric. \qed
\medskip

{\bf Proposition 5.9.} {\sl Let $(X,*)$ be a group with identity $e$.  Then every subgroup $(A,*)$ of $(X,*)$ is divisible.}
\medskip

{\it Proof.} Let $x, y\in A$ such that $a*x = y$ for some $a\in X$.  It follows that $a= y*x^{-1}\in A$ since $(A,*)$ is a subgroup of $X$.  This shows that $(A,*)$ is divisible in $(X,*)$. \qed

\bigskip

\footnotesize{\textsc{H. Fayoumi, Department of Mathematics and Statistics,
University of Toledo, Toledo, OH 43606-3390, U. S. A.}
\par \textit{E-mail address}: \texttt{hiba.fayoumi@UToledo.edu}\par
\medskip

 \footnotesize{\textsc{ Hee Sik Kim, Department of Mathematics, Hanyang University, Seoul, 04763, Korea}
\par \textit{E-mail address}: \texttt{heekim@hanyang.ac.kr}}
\medskip

\bigskip
%%%%%%%%%%%%%%%%%%%%%%%%%%%%%%%%%%%%%%%%%%%%%%%%%%%%%%%%%%%%%

\end{document}